\documentclass[11pt]{article}
\usepackage{amssymb,authblk}
\usepackage{amsmath}
\usepackage{amsthm}
\usepackage{bm}
\usepackage{xcolor}

\usepackage{geometry}
\usepackage[colorlinks=true,
linkcolor=blue,citecolor=blue,
urlcolor=blue]{hyperref}

\newtheorem{lemma}{Lemma}[section]
\newtheorem{theorem}{Theorem}[section]

\numberwithin{equation}{section}
\allowdisplaybreaks

\begin{document}
\title{A characterization of extremal non-transmission-regular graphs by the distance (signless Laplacian) spectral radius}
\author[a]{Jingfen Lan\thanks{Corresponding author}}
\author[b]{Lele Liu}
\affil[a]{School of Mathematics and Statistics, Xidian University, Xi'an 710126, China. (\texttt{jflan@xidian.edu.cn})}
\affil[b]{School of Mathematical Sciences, Anhui University, Hefei 230601, China. (\texttt{liu@ahu.edu.cn})}
\date{}
\maketitle

\begin{abstract}
Let $G$ be a simple connected graph of order $n$ and $\partial(G)$ is the spectral radius
of the distance matrix $D(G)$ of $G$. The transmission $D_i$ of vertex $i$ is the $i$-th row
sum of $D(G)$. Denote by $D_{\max}(G)$ the maximum of transmissions over all vertices of $G$,
and $\partial^Q(G)$ is the spectral radius of the distance signless Laplacian matrix $D(G)+\mbox{diag}(D_1,D_2,\ldots,D_n)$.
In this paper, we present a sharp lower bound of $2D_{\max}(G)-\partial^Q(G)$ among all $n$-vertex
connected graphs, and characterize the extremal graphs. Furthermore, we give the minimum
values of respective $D_{\max}(G)-\partial(G)$ and $2D_{\max}(G)-\partial^Q(G)$ on trees
and characterize the extremal trees.

\end{abstract}

\section{Introduction}
\label{intro}
All graphs considered in this paper are undirected and simple. Let $G$ be a connected graph
of order $n$ with vertex set $V(G)=\{v_1,v_2,\ldots,v_n\}$. The distance between two vertices
$v_i$ and $v_j$, denoted by $d_{ij}$, is the length of the shortest path between $v_i$ and $v_j$.
The distance matrix $D(G)$ of graph $G$ is the $n\times n$ matrix with $d_{ij}$ as the $(i, j)$-entry,
whose largest eigenvalue, denoted by $\partial(G)$, is called the \emph{distance spectral radius}
of $G$. The \emph{transmission} $D_i(G)$ of vertex $v_i$ is the sum of distance from $v_i$ to all
over vertices of $G$, i.e., the $i$-th row sum of $D(G)$. Denote by $D_{\max}(G)$ and $D_{\min}(G)$
respectively the maximum and minimum of the transmissions over all vertices of $G$. If there is
no risk of confusion, we denote $D_i(G)$, $D_{\max}(G)$ and $D_{\min}(G)$ simply by $D_i$, $D_{\max}$
and $D_{\min}$ respectively. Let $Q(G)=D(G)+\mbox{diag}(D_1, D_2,\ldots, D_n)$ be the
\emph{distance signless Laplacian matrix} of graph $G$, whose largest eigenvalue, denoted by $\partial^Q(G)$,
is called the \emph{distance signless Laplacian  spectral radius} of $G$. The
\emph{distinguished vertex deleted regular graph} (DVDR), introduced by Atic and Panigrahi \cite{AP},
is a connected graph $G$ with a vertex $v$ such that the degree $d(v)$ of $v$ equals to $n-1$ and
$G-v$ is regular. The vertex $v$ is called a \emph{distinguished vertex} of the DVDR graph $G$.
If $G-v$ is $r$-regular, we say $G$ is an $r$-DVDR graph. For example, wheel graphs are 2-DVDR graphs
and complete graphs are $(n-1)$-DVDR graphs. The \emph{Wiener index} $W(G)$ of a graph $G$ is the
sum of distances between all unordered pairs of vertices of $G$, that is, the half sum of all the entries of $D(G)$,
\[
W(G) = \frac{1}{2} \sum_{i=1}^n D_i.
\]

It is well-known that the adjacency spectral radius, denoted by $\rho(G)$, is no more than the
maximum degree $\Delta(G)$ of a graph $G$, and $\rho(G)=\Delta(G)$ if and only if $G$ is
regular for a connected graph $G$ (see \cite{CS}). Hence, $\Delta(G)-\rho(G)$ can be
regarded as a measure of irregularity for a nonregular graph $G$ (\cite{CR}, p.\,242),
and many estimates on it have been obtained (\cite{CGN,LHY,Z,S,CH}). Similarly,
$2\Delta(G) - q(G)$ is also used to qualify how close a graph to being regular (\cite{NLL,SPS}),
where $q(G)$ is the signless Laplacian spectral radius of $G$. As a natural generalization,
graph regularity associated with the distance matrix is defined. A connected graph $G$ is
\emph{transmission-regular} if $D_{\max} = D_{\min}$; otherwise $G$ is \emph{non-transmission-regular}.
It is known that the distance (signless Laplacian) spectral radius of a graph $G$ satisfies
$D_{\min}(G)\leq \partial(G)\leq D_{\max}(G)$ (and $2D_{\min}(G)\leq \partial^Q(G)\leq 2D_{\max}(G)$),
with the equalities if and only if $G$ is a transmission-regular graph \cite{M}.
Naturally, Liu, Shu and Xue \cite{LSX} asked the question:
how small can $D_{\max}(G)-\partial (G)$ and $2D_{\max}(G)-\partial^Q(G)$ be for non-transmission-regular
graphs? They gave some lower bounds for the two indicators. In \cite{LHC}, Liu, Shan and He confirmed
a conjecture with respect to $D_{\max}(G)-\partial(G)$ posed by Liu, Shu and Xue \cite{LSX}
by giving a tight lower bound for graphs $G$ of order $n$, and determining the extremal
graphs to be DVDR graphs. For a connected graph $G$ of order $n$, write
\begin{align*}
\sigma(G) & := D_{\max}(G) - \partial(G), \\
\tau(G) & := 2D_{\max}(G) - \partial^Q(G).
\end{align*}
Using the method in \cite{LHC}, we give a tight lower bound for $\tau(G)$ and characterize the
extremal graphs which are the same DVDR graphs as in the case of $\sigma(G)$. Furthermore,
we give lower bounds respectively for $\sigma(G)$ and $\tau(G)$ on trees and characterize the
extremal graphs which are stars.

\begin{theorem}\label{mainthm1}
Let $G$ be a connected non-transmission-regular graph on $n$ vertices.
\begin{enumerate}
\item[$(1)$] If $n$ is odd, then
\[
\tau(G)\geq \frac{n+2-\sqrt{n^2+4n-4}}{2}.
\]
Equality holds if and only if $G\cong K_{1,2,\ldots,2}$.

\item[$(2)$] If $n$ is odd, then
\[
\tau(G)\geq \frac{n+4-\sqrt{n^2+8n}}{2}.
\]
Equality holds if and only if $G$ is isomorphic to some $(n-4)$-DVDR graph.
\end{enumerate}
\end{theorem}

Consider $\sigma(G)$ and $\tau(G)$ only on trees, which are also non-transmission-regular graphs.
We obtain the lower bounds and the extremal trees.

\begin{theorem}\label{mainthm2}
Let $T$ be a tree on $n$ vertices. We have
\[
\sigma(T)\geq \frac{2n-2-\sqrt{4n^2-12n+12}}{2}
\]
and
\[
\tau(T)\geq \frac{3n-4-\sqrt{9n^2-32n+32}}{2},
\]
with the equalities if and only if $T\cong K_{1,n-1}$.
\end{theorem}

\section{Preliminaries}
\label{sec:1}
Let $M$ be a real symmetric matrix associated with a graph $G$. An equitable partition
$\pi : V(G)=X_1\cup X_2\cup \cdots \cup X_k$ of the vertex set $V(G)$ can divide $M$ into
the block form $M=(M_{ij})$ such that any block $M_{ij}$ has a constant row sum $b_{ij}$,
and then the eigenvalues of the $k\times k$ matrix $B=(b_{ij})$ are also eigenvalues of $M$ (see \cite[p.\,24]{Brouwer2011}).
Furthermore, if $M$ is nonnegative and irreducible then $B$ and $M$ have the same spectral radius. Based
on this understanding, let us calculate $\tau(K_{1,2,\ldots,2})$ and $\tau(G)$ of an $(n-4)$-DVDR graph $G$.

\begin{lemma}\label{lm2.1}
Let $K_{1,2,\ldots,2}$ be the $r$-partite graph on $n$ vertices with $n=2r-1$ and $G$ is
an $(n-4)$-DVDR graph. Then
\[
\tau(K_{1,2,\ldots,2}) = \frac{n+2-\sqrt{n^2+4n-4}}{2}~\ \text{and}\ \tau(G)=\frac{n+4-\sqrt{n^2+8n}}{2}.
\]
\end{lemma}

\begin{proof}
Note that $K_{1,2,\ldots,2}$ is an $(n-3)$-DVDR graph. Denote by $v$ the distinguish vertex
of $K_{1,2,\ldots,2}$ (or $G$). The quotient matrix of $Q(K_{1,2,\ldots,2})$ with respect
to the equitable partition $\{v\}\cup (V(K_{1,2,\ldots,2})\backslash \{v\})$ is
\[
\begin{bmatrix}
  n-1 & n-1 \\ 1 & 2n-1
\end{bmatrix},
\]
which has the same spectral radius as $Q(K_{1,2,\ldots,2})$, i.e.,
\[
\partial^Q(K_{1,2,\ldots,2}) = \frac{3n-2+\sqrt{n^2+4n-4}}{2}.
\]
Since $D_{\max}(K_{1,2,\ldots,2})=n$, we have
\[
\tau(K_{1,2,\ldots,2}) = 2n-\partial^Q(K_{1,2,\ldots,2})
= \frac{n+2-\sqrt{n^2+4n-4}}{2}.
\]

The quotient matrix of $Q(G)$ with respect to the equitable partition $\{v\}\cup (V(G)\backslash \{v\})$ is
\[
\begin{bmatrix}
  n-1 & n-1 \\ 1 & 2n+1
\end{bmatrix},
\]
which has the same spectral radius as $Q(G)$, i.e.,
\[
\partial^Q(G) = \frac{3n+\sqrt{n^2+8n}}{2}.
\]
Since $D_{\max}(G) = n+1$, we have
\[
\tau(K_{1,2,\ldots,2}) = 2(n+1)-\partial^Q(G)
= \frac{n+4-\sqrt{n^2+8n}}{2}.
\]
The proof is completed.
\end{proof}

In the following, we denote
\[
\tau_n:= \frac{n+2\gamma_n-\sqrt{(n+2\gamma_n)^2-8\gamma_n}}{2},
\]
where $\gamma_n=1$ if $n$ is odd; $\gamma_n=2$ if $n$ is even. We can see that $\gamma_n$ satisfies the equality
\begin{equation}\label{eq1}
\tau_n^2 - (n+2\gamma_n)\tau_n + 2\gamma_n = 0.
\end{equation}

By a similar calculation as the proof of Lemma \ref{lm2.1}, we obtain $\sigma(K_{1,n-1})$ and $\tau(K_{1,n-1})$.

\begin{lemma}\label{lm2.2}
Let $K_{1,n-1}$ be the star graph on $n$ vertices. We have
\[
\sigma(K_{1,n-1})=\frac{2n-2-\sqrt{4n^2-12n+12}}{2} \ \text{and}\ \tau(K_{1,n-1})=\frac{3n-4-\sqrt{9n^2-32n+32}}{2}.
\]
\end{lemma}

We denote
\[
\sigma'_n:= \frac{n+\eta_n-\sqrt{(n+\eta_n)^2-4\eta_n}}{2},
\]
\[
\tau'_n:= \frac{n+2\eta_n-\sqrt{(n+2\eta_n)^2-8\eta_n}}{2}
\]
with $\eta_n=n-2$. We see that $\sigma'_n$ and $\tau'_n$ satisfy the following equalities
\begin{align}\label{eq2}
\sigma'^2_n -(n+\eta_n)\sigma'_n+\eta_n & = 0, \\
\tau'^2_n -(n+2\eta_n)\tau'_n+2\eta_n & = 0. \label{eq3}
\end{align}

\section{Proof of Theorem \ref{mainthm1}}
\label{sec:2}
In this section, we assume that $G^*$ is a graph attaining the minimum of $\tau(G)$ among
connected non-transmission-regular graphs $G$ of order $n$. Let $\bm{x}$ be the unit positive
eigenvector of $Q(G^*)$ corresponding to $\partial^Q(G^*)$. For convenience, we denote
$x_{\max}:= \max \{x_v\,|\,v\in V(G^*)\}$ and $x_{\min}:= \min \{x_v\,|\,v\in V(G^*)\}$. Let
$u$ and $v$ be respectively the vertices such that $x_u = x_{\max}$ and $x_{v} = x_{\min}$.
By the equation $Q(G^*) \bm{x} = \partial^Q(G^*) \bm{x}$ with respect to the vertex $u$, we see
\begin{align*}
\partial^Q(G^*) x_u & = D_ux_u+\sum_{w\in V(G^*)} d(u, w)x_w \\
& = D_ux_u+\sum_{w\in V(G^*)\backslash\{v\}} d(u,w)x_w + d(u,v)x_v \\
&\leq (2D_{\max} - d(u,v))x_u + d(u,v)x_v,
\end{align*}
which yields that
\[
2D_{\max} - \partial^Q(G^*) \geq d(u,v)(1-\frac{x_v}{x_u}) \geq 1-\frac{x_v}{x_u}.
\]
On the other hand, $\tau(G^*) = 2D_{\max} - \partial^Q(G^*) \leq \tau_n$, we obtain
\begin{equation}\label{eq4}
\frac{x_{\max}}{x_{\min}} = \frac{x_u}{x_v}\leq \frac{1}{1-\tau_n}.
\end{equation}

Theorem~\ref{mainthm1} follows from Lemmas \ref{lm3.1} and \ref{lm3.2} below.

\begin{lemma}\label{lm3.1}
$\tau(G^*) = \tau_n$.
\end{lemma}

\begin{proof}
By the minimality of $\tau(G^*)$, we have $\tau(G^*)\leq\tau_n$. Noting that
\[
(2D_{\max}-\partial^Q(G^*)) \sum_{i=1}^n x_i = 2\sum_{i=1}^n (D_{\max}-D_i)x_i \geq 2(nD_{\max}-2W) x_{\min},
\]
together with Eqs. \eqref{eq1} and \eqref{eq4}, we find
\begin{equation}\label{eq5}
\begin{split}
nD_{\max} - 2W & \leq \frac{\tau_n(x_{\min}+(n-1) x_{\max})}{2x_{\min}} \\
& = \frac{\tau_n}{2} \Big(1 + (n-1) \cdot \frac{x_{\max}}{x_{\min}}\Big) \\
& \leq \frac{\tau_n(n-\tau_n)}{2(1-\tau_n)} = \gamma_n.
\end{split}
\end{equation}
Hence, $nD_{\max} - 2W = 1$ if $n$ is odd; If $n$ is even, then $nD_{\max}-2W$ is an even
positive integer and equals to $2$. Anyway, we have
\begin{equation}\label{eq6}
nD_{\max} - 2W = \gamma_n.
\end{equation}
Thus, the first equality of Eq. \eqref{eq5} holds and implies $\tau(G^*)=\tau_n$.
\end{proof}

\begin{lemma}\label{lm3.2}
If $n$ is odd, then $G^*\cong K_{1,2,\ldots,2}$; If $n$ is even, then $G^*$ is isomorphic
to some $(n-4)$-DVDR graph.
\end{lemma}

\begin{proof}
By the first equality of Eq. \eqref{eq5}, there are $n-1$ vertices attaining $x_{\max}$ and
the remaining vertex attains $x_{\min}$. Our proof relies on the the following two claims.

\noindent{\bf Claim 1.} There are $n-1$ vertices attaining $D_{\max}$.

Let $S=\{w\in V(G^*): D_w = D_{\max}\}$. From Eq. \eqref{eq6}, we infer that $|S|=n-1$
if $n$ is odd, and $|S|\in \{n-1, n-2\}$ if $n$ is even. Suppose that $|S|=n-2$, in
which case the remaining two vertices both attain $D_{\min}$. There must be a vertex $u$
such that $x_u = x_{\max}$ and $D_u = D_{\min}$. However, the eigenvalue equation
\[
2D_{\min} x_u < \partial^Q(G^*) x_u = \sum_{w\in V(G^*)} d(u,w)x_w + D_ux_u \leq 2D_ux_u
\]
implies $D_{\min} < D_u$. A contradiction! Thus, $|S|=n-1$ for either odd or even $n$.

\noindent{\bf Claim 2.} $D_{\max} = n + \gamma_n - 1$ and $D_{\min} = n-1$.

By the proof of Claim 1, the $n-1$ vertices attaining $x_{\max}$ also attain $D_{\max}$.
Recalling that $x_v=x_{\min}$, we have $D_v = D_{\min} = D_{\max}-\gamma_n$ by Eq. \eqref{eq6}.
It can be derived from the eigenvalue equation that
\begin{align*}
(2D_{\max} - \tau_n) x_{\min}
& = \partial^Q(G^*) x_v = \sum_{u\in V(G^*)} d(v,u)x_u + D_vx_v \\
& = D_vx_u + D_vx_v = (D_{\max}-\gamma_n)(x_{\max} + x_{\min}).
\end{align*}
Besides, the equalities in Eq. \eqref{eq5} imply $\frac{x_{\max}}{x_{\min}} = \frac{1}{1-\tau_n}$.
Together with Eq. \eqref{eq1}, it follows that
\begin{align*}
D_{\max} & = \frac{\gamma_n(1+\frac{x_{\max}}{x_{\min}})-\tau_n}{\frac{x_{\max}}{x_{\min}}-1}
= \frac{\frac{\gamma_n(2-\tau_n)}{1-\tau_n}-\tau_n}{\frac{\tau_n}{1-\tau_n}} \\
& = \tau_n-\gamma_n+\frac{2\gamma_n}{\tau_n}-1 = n-1+\gamma_n.
\end{align*}
Hence, $D_{\min} = D_{\max} - \gamma_n = n-1$. \par\vspace{2mm}

Now we continue the proof of the lemma. Since $D_v=D_{\min}=n-1$, we have $d(v)=n-1$,
$d(u,v)=1$ and $d(u,w)\leq 2$ for all $u,w \in V(G^*)\backslash \{v\}$, which yields that
\[
D_u = \sum_{w\in V(G^*)} d(u,w) = d(u) + 2(n-1-d(u)) = 2(n-1) - d(u).
\]
Therefore, for each vertex $u \in V(G^*)\backslash \{v\}$,
\[
d(u) = 2(n-1) - D_u = n - 1 - \gamma_n =
\begin{cases}
n-2, & n\ \text{is odd;} \\ n-3, & n\ \text{is even.}
\end{cases}
\]
Hence, $G^*\cong K_{1,2,\ldots,2}$ if $n$ is odd, and $G^*$ is isomorphic to some
$(n-4)$-DVDR graph if $n$ is even.
\end{proof}

\section{Proof of Theorem \ref{mainthm2}}
\label{sec:3}
Suppose that $D_{\max}$ and $W$ are respectively the maximum transmission and
Wiener index of a tree $T$. It is easy to see that
\[
nD_{\max} - 2W = \sum_{i=1}^n (D_{\max}-D_i).
\]
Let $e=v_1v_2$ be a leaf of $T$ with $d(v_1)=1$. Since $d(v_1,v)=d(v_2,v)+1$ for any
vertex $v\in V(T)\backslash \{v_1, v_2\}$, we have $D_{v_1}-D_{v_2}=n-2$. Thus,
\begin{equation}\label{eq7}
nD_{\max} - 2W \geq D_{\max} - D_{v_2} \geq D_{v_1} - D_{v_2}
= n-2 = \eta_n.
\end{equation}

We prove Theorem \ref{mainthm2} for cases of $\tau(T)$ and $\sigma(T)$ separately
by Lemmas \ref{lm4.1}~and \ref{lm4.2} below.

Assume that $T^*$ is a tree attaining the minimum of $\tau(T)$ among all trees $T$ of
order $n$. Let $\bm{x}$ be the unit positive eigenvector of $Q(T^*)$ corresponding to $\partial^Q(T^*)$.
For convenience, we still denote $x_{\max}:=\max \{x_v\,|\,v\in V(T^*)\}$ and $x_{\min}:=\min \{x_v\,|\,v\in V(T^*)\}$.
Let $u$ and $v$ be respectively the vertices such that $x_u = x_{\max}$ and $x_{v} = x_{\min}$.
By a similar derivation of Eq. \eqref{eq4}, we get
\begin{equation}\label{eq8}
\frac{x_{\max}}{x_{\min}} = \frac{x_u}{x_v}\leq
\frac{1}{1-\tau'_n}.
\end{equation}

\begin{lemma}\label{lm4.1}
$\tau(T^*) = \tau'_n$ and $T^*\cong K_{1,n-1}$.
\end{lemma}

\begin{proof}
By the minimality of $\tau(T^*)$, we have $\tau(T^*)\leq\tau'_n$. Noting that
\[
(2D_{\max} - \partial^Q(T^*)) \sum_{i=1}^n x_i = 2\sum_{i=1}^n (D_{\max} - D_i)x_i \geq  2(nD_{\max} - 2W)x_{\min},
\]
together with Eqs. \eqref{eq3} and \eqref{eq8}, we see
\begin{equation}\label{eq9}
\begin{split}
nD_{\max} - 2W
& \leq \frac{\tau_n'(x_{\min} + (n-1) x_{\max})}{2x_{\min}} \\
& = \frac{\tau'_n}{2} \Big(1 + (n-1)\cdot\frac{x_{\max}}{x_{\min}}\Big) \\
& \leq \frac{\tau'_n(n-\tau'_n)}{2(1-\tau'_n)} = \eta_n.
\end{split}
\end{equation}
Hence, by Eq.~\eqref{eq7}, we have
\begin{equation}\label{eq10}
nD_{\max} - 2W = \eta_n.
\end{equation}
Thus, all equalities of Eq. \eqref{eq9} hold. We get $\tau(T^*)=\tau_n$ and there are $n-1$ vertices attaining
$x_{\max}$ and the remaining vertex attains $x_{\min}$ in $T^*$. For the same reason of Claim 1 of Lemma \ref{lm3.2},
the $n-1$ vertices attaining $x_{\max}$ do not attain $D_{\min}$. By the argument of Eq. \eqref{eq7}, the
equality \eqref{eq10} implies that $v$ is the unique vertex attaining $D_{\min}$ and the other $n-1$ vertices
attain $D_{\max}$, i.e., $D_v = D_{\min} = D_{\max} - \eta_n$. Noting the similarity of Eq. \eqref{eq3}
and Eq. \eqref{eq1}, a similar deduction as Claim~2 of Lemma \ref{lm3.2} can conclude
$D_{\max} = n-1+\eta_n$ and $D_{\min} = D_{\max} - \eta_n = n-1$.

For the tree $T^*$, $D_v=D_{\min}=n-1$ forces $d(u,v)=1$, which means $v$ is adjacent to any other vertex
of $T^*$. Hence, $T^*\cong K_{1,n-1}$.
\end{proof}

In the following, we still denote by $T^*$ an extremal tree attaining the minimum of $\sigma(T)$ among
all trees of order $n$ and $\bm{x}$ is the unit positive eigenvector of $D(T^*)$ corresponding to
$\partial(T^*)$. Let the vertices $u, v$ such that $x_u = x_{\max}$ and $x_{v} = x_{\min}$ have the
same meaning as before. By the equation $D(G^*) \bm{x} = \partial \bm{x}$ with respect to the vertex
$u$, we see
\begin{align*}
\partial(T^*) x_u & = \sum_{w\in V(G^*)} d(u, w)x_w \\
& = \sum_{w\in V(G^*)\backslash\{v\}} d(u,w)x_w + d(u,v)x_v \\
& \leq (D_{\max} - d(u,v)) x_u + d(u,v)x_v,
\end{align*}
which yields that
\[
D_{\max}-\partial(T^*) \geq d(u,v) \Big(1 - \frac{x_v}{x_u} \Big) \geq 1-\frac{x_v}{x_u}.
\]
On the other hand, $\sigma(T^*) = D_{\max} - \partial(T^*) \leq \sigma'_n$, we obtain
\begin{equation}\label{eq11}
\frac{x_{\max}}{x_{\min}} = \frac{x_u}{x_v}\leq \frac{1}{1 - \sigma'_n}.
\end{equation}

\begin{lemma}\label{lm4.2}
$\sigma(T^*) = \sigma'_n$ and $T^*\cong K_{1,n-1}$.
\end{lemma}

\begin{proof}
By the eigenvalue equation
\[
(D_{\max} - \partial(T^*))\sum_{i=1}^n x_i = \sum_{i=1}^n (D_{\max}-D_i) x_i \geq (nD_{\max} - 2W) x_{\min},
\]
together with Eqs. \eqref{eq2} and \eqref{eq10}, we find
\begin{equation}\label{eq12}
\begin{split}
nD_{\max} - 2W
& \leq \frac{\sigma'_n(x_{\min} + (n-1)x_{\max})}{x_{\min}} \\
& = \sigma'_n \Big(1 + (n-1) \cdot \frac{x_{\max}}{x_{\min}}\Big) \\
& \leq \frac{\sigma'_n(n-\sigma'_n)}{1-\sigma'_n}) = \eta_n.
\end{split}
\end{equation}
Hence, by Eq. \eqref{eq7}, we have
\begin{equation}\label{eq13}
nD_{\max} - 2W = \eta_n.
\end{equation}
Thus, all equalities of Eq. \eqref{eq12} hold. We get $\sigma(T^*)=\sigma'_n$ and there are $n-1$
vertices attaining $x_{\max}$ and the remaining vertex attains $x_{\min}$ in $T^*$. By the argument
of Eq. \eqref{eq7}, the equality \eqref{eq13} implies that $v$ is the unique vertex attaining $D_{\min}$
and the other $n-1$ vertices attain $D_{\max}$, i.e., $D_v = D_{\min} = D_{\max} - \eta_n$. Using
the eigenvalue equation, we have
\[
(D_{\max} - \sigma'_n) x_{\min} = \partial(T^*) x_v = \sum_{u\in V(G^*)} d(v,u)x_u
= D_vx_u = (D_{\max} - \eta_n) x_{\max}.
\]
The equalities in Eq. \eqref{eq12} imply $\frac{x_{\max}}{x_{\min}} = \frac{1}{1-\sigma'_n}$.
Together with Eq. \eqref{eq2}, it follows that
\[
D_{\max} = \frac{\eta_n\cdot\frac{x_{\max}}{x_{\min}} - \sigma'_n}{\frac{x_{\max}}{x_{\min}} - 1}
= \frac{\frac{\eta_n}{1-\sigma'_n}-\sigma'_n}{\frac{\sigma'_n}{1-\sigma'_n}}
= \sigma'_n + \frac{\eta_n}{\sigma'_n}-1 = n - 1 + \eta_n.
\]
Hence, $D_v = D_{\min} = D_{\max} - \eta_n = n-1$, which means $T^*\cong K_{1,n-1}$.
\end{proof}

\section{Concluding Remark}
\label{sec:4}
Although it is easy to see that $\sigma(G)<\tau(G)$ for non-transmission-regular graphs,
the results of this paper show that the extremal graphs attaining minimum $\sigma(G)$ and
those attaining minimum $\tau(G)$ are the same, and the extremal trees attaining minimal
$\sigma(T)$ and those attaining minimal $\tau(T)$ are also the same. This may mean that
the two indicators for measuring non-transmission-regular graphs are similar to some
extent. In addition, one can find that the the minimum values $\sigma_n$ \cite{LHC} and
$\tau_n$ decrease and converge to zero as $n$ increases to infinity, while $\sigma'_n$
and $\tau'_n$ increase as $n$ increases. This means that there is a big gap between
trees and cyclic graphs in terms of $\sigma$ and $\tau$. However, we know that the
stars take the maximal value of $\Delta(G)-\rho(G)$. Our results reveal that $\sigma$
and $\tau$ are quite different from $\Delta(G)-\rho(G)$ in measuring the irregularity
of graphs.


%
%



\end{document}